\documentclass[12pt,a4paper]{amsart}

\def\marginpar#1{}
\def\resizebox#1#2#3{}

\newcommand {\supplus}{\mathop{{\supset}\llap{\raise
0.5pt\hbox{\normalfont\small+}\hskip 0.5pt}}}

\newcommand {\subplus}{\mathop{{\subset}\llap{\raise
0.5pt\hbox{\normalfont\small+}\hskip 0.5pt}}}

\newcommand {\Cee}    {{\mathbb  C}}

\newcommand {\Ree}    {{\mathbb  R}}

\newcommand {\Zee}    {{\mathbb  Z}}

\newcommand {\fder}   {{\mathfrak{der}}}   %

\newcommand {\fg}     {{\mathfrak{g}}}    %
\newcommand {\fgl}    {{\mathfrak{gl}}}  %
\newcommand {\fh}     {{\mathfrak{h}}}

\newcommand {\fl}     {{\mathfrak{l}}}
\newcommand {\fL}     {{\mathfrak{L}}}

\newcommand {\fn}     {{\mathfrak{n}}}

\newcommand {\fo}     {{\mathfrak{o}}}

\newcommand {\fp}    {{\mathfrak{p}}}   %

\newcommand {\fsl}    {{\mathfrak{sl}}}

\newcommand {\fsp}    {{\mathfrak{sp}}}

\newcommand {\fvect}  {{\mathfrak{vect}}}   %

\newcommand {\cal} {\mathcal}

\newcommand {\cF}     {{\cal F}}

\def \opname#1#2%
  {\expandafter\newcommand \csname #1\endcsname {{\mathop{#2}\nolimits}}}

\newcommand{\rmname}[1]
  {\expandafter\newcommand \csname #1\endcsname {{\operatorname{#1}}}}

\newcommand{\rmnameii}[2]
  {\expandafter\newcommand \csname #1\endcsname {{\operatorname{#2}}}}

\rmname{act}
\rmname{Ad}
\rmname{Add}
\rmname{ad}
\rmname{Alt}
\rmname{alt}
\rmname{Ann}
\rmname{antidiag}
\rmname{Ber}
\rmname{ber}
\rmname{Br}
\rmname{card}
\rmname{ch}
\rmname{Char}
\rmname{cem}
\rmname{cj}
\rmname{Cliff}
\rmname{cntr}
\rmname{codim}
\rmname{coind}
\rmname{const}
\rmname{col}
\rmname{cork}
\rmname{cpr}
\rmname{diag}
\rmnameii{Div}{div}
\rmname{Def}
\rmname{Der}
\rmname{Dim}
\rmname{End}
\rmname{Even}
\rmname{Ext}
\rmname{gr}
\rmname{Hom}
\rmname{HT}
\rmnameii{Ht}{ht}
\rmname{hwt}
\rmname{Id}
\rmname{id}
\rmname{ind}
\rmname{Ind}
\rmname{Coind}
\rmname{Inf}
\rmname{irr}
\rmname{Le}
\rmname{Lie}
\rmname{lwt}
\rmname{mult}
\rmname{Mor}
\rmname{nm}
\rmname{Ob}
\rmname{Odd}
\rmname{Osc}
\rmname{per}
\rmname{Pic}
\rmname{pr}
\rmname{pro}
\rmname{Prime}
\rmname{Proj}
\rmname{prt}
\rmname{pt}
\rmname{Q}
\rmname{qet}
\rmname{qtr}
\rmname{rd}
\rmname{rk}
\rmname{row}
\rmname{Res}
\rmname{salt}
\rmname{Sch}
\rmname{SBr}
\rmname{scalar}
\rmname{Ser}
\rmname{sign}
\rmname{Smbl}
\rmname{spin}
\rmname{ssym}
\rmname{str}
\rmname{st}
\rmname{sgn}
\rmname{sq}
\rmname{symm}
\rmname{supp}
\rmname{Supp}
\rmname{St}
\rmname{Spec}
\rmname{Spm}
\rmname{tr}
\rmname{vpt}
\rmname{weyl}
\rmname{Weyl}
\rmname{Witt}

\opname{vvol}  {{v\hspace{-0.1ex}o\hspace{-0.02ex}l\/}}
\opname{pnt}  {\text{\normalfont pt}}
\rmname{Span} 
\opname{slim} {\overline{\lim}}
\opname{Vol}  {{V\hspace{-0.55ex}o\hspace{-0.02ex}l\/}}
\opname{QVol} {{Q\hspace{-0.3ex}V\hspace{-0.55ex}o\hspace{-0.02ex}l\/}}
\opname{PoVol}{{P\hspace{-0.35ex}o\hspace{-0.25ex}V\hspace{-0.55ex}o\hspace{-0.02ex}l\/}}
\opname{BVol} {{B\hspace{-0.2ex}V\hspace{-0.55ex}o\hspace{-0.02ex}l\/}}
\opname{Par}  {{P\hspace{-0.3ex}a\hspace{-0.05ex}r\/}}

\rmname{Mat}
\rmname{Bil}
\rmname{Diff}
\rmname{Ker}
\rmname{Herm}
\rmname{Coker}
\rmname{Conn}
\rmname{Covect}
\rmname{Vect}
\rmname{Int}

\rmnameii {IM} {Im}
\rmnameii {RE} {Re}

\opname{Aut} {{A\hspace{-0.2ex}u\hspace{-0.1ex}t\/}}
\opname{GL} {{G\hspace{-0.3ex}L}}
\opname{SL} {{S\hspace{-0.3ex}L}}
\opname{Exp} {{E\hspace{-0.2ex}x\hspace{-0.1ex}p\/}}
\opname{GQ} {{G\hspace{-0.2ex}Q}}
\opname{OSp} {{O\hspace{-0.25ex}S\hspace{-0.15ex}p\/}}
\opname{Out} {{O\hspace{-0.25ex}u\hspace{-0.15ex}t\/}}
\opname{Spp} {{S\hspace{-0.2ex}p\/}}
\opname{SpO} {{S\hspace{-0.2ex}p\hspace{-0.02ex}O\/}}
\opname{Pe} {{P\hspace{-0.25ex}e\/}}
\opname{SPe} {{S\hspace{-0.25ex}P\hspace{-0.25ex}e\/}}
\opname{Spin} {{S\hspace{-0.25ex}p\hspace{-0.05ex}i\hspace{-0.1ex}n\/}}
\opname{Iso} {{I\hspace{-0.25ex}s\hspace{-0.1ex}o\/}}
\opname{SSPe} {{S\hspace{-0.25ex}S\hspace{-0.15ex}P\hspace{-0.25ex}e\/}}
\opname{PeU} {{P\hspace{-0.25ex}e\hspace{-0.1ex}U\/}}
\opname{QU} {{Q\hspace{-0.15ex}U\/}}
\opname{U} {{U\/}}

\opname{cGQ} {{\cal G \hspace{-0.2em} Q \/}}
\opname{cSL} {{\cal S \hspace{-0.2em} L \/}}
\opname{cGL} {{\cal G \hspace{-0.2em} L \/}}\opname{cGr} {{\cal G \hspace{-0.2em} r \/}}
\opname{cOSp} {{\cal O \hspace{-0.2em} S \hspace{-0.3em} \it p\/}}
\opname{cPe} {{\cal P \hspace{-1.5pt} \it e\/}}
\opname{cVect} {{\cal V \hspace{-1.5pt} \it e\hspace{-0.1ex}c\hspace{-0.1ex}t\/}}
\opname{cVol} {{\cal V \hspace{-1.5pt} \it o\hspace{-0.1ex}l\/}}
\opname{cAut} {{\cal A \hspace{-0.2em} \it u\hspace{-0.1em}t\/}}
\opname{cCovect} {{\cal C \hspace{-1.5pt}
     \it o\hspace{-0.1ex}v\hspace{-0.1ex}e\hspace{-0.1ex}c\hspace{-0.1ex}t\/}}
\opname{CW} {{C\hspace{-0.15ex}W}}

\newcommand {\tto} {\longrightarrow}
\newcommand {\pder}[1] {{\frac{\partial}{\partial {#1}}}}

\newcommand {\bcdot}   {\mathbin{\hbox{\raise.4ex\hbox{\bf.}}}} 

\newcommand {\secno} {}
\newcommand {\ssecfont} {\normalfont\bf}

\newtheorem*{Theorem}{\secno Theorem}

\newtheorem*{Lemma}{\secno Lemma}

\newenvironment {th*}[1]
    {\gdef\thname{#1} \begin{thn}}%
    {\end{thn}}
\newtheorem*{thn} {\thname}

\theoremstyle{definition}

\newenvironment {ex*}[1]
    {\gdef\thname{#1} \begin{exn}}%
    {\end{exn}}
\newtheorem*{exn}{\thname}

\theoremstyle{remark}

\newenvironment {rem*}[1]
    {\gdef\thname{#1} \begin{remn}}%
    {\end{remn}}
\newtheorem*{remn}{\thname}

\newcommand {\ssec}{\subsection*}

\newcommand {\ssbegin}[2]
  {\def \secno {\gdef \secno {}{\ssecfont #1. }}%
   \begin{#2}}

\begin{document}

\title[Lie superalgebras in cohomology]{Lie superalgebra structures in
$C^{\bcdot}(\fg;\fg)$ and $H^{\bcdot}(\fg;\fg)$}

\author{Alexei Lebedev$^1$, Dimitry Leites$^2$, Ilya Shereshevskii$^1$}

\address{$^1$Institute for Physics of Microstructures, RAS, GSP-105, RU-603950,
Nizhny Novgorod, Russia; yorool@mail.ru;ilya@ipm.sci-nnov.ru;
$^2$Department
of Mathematics, University of Stockholm, Roslagsv.  101, Kr\"aftriket
hus 6, SE-106 91, Stockholm, Sweden; mleites@ma\-te\-ma\-tik.su.se}

\thanks{For financial support and stimulating working conditions D.L.
is thankful to MPIM-Bonn and IHES}



\date{October 31, 2003}

\dedicatory{To Arkady Lvovich Onishchik.}

\begin{abstract} Let $\fn=\fvect(M)$ be the Lie (super)algebra of
vector fields on any connected (super)manifold $M$; let $\Pi$ be the
change of parity functor, $C^{i}$ and $H^{i}$ the space of $i$-chains
and $i$-cohomology.  The Nijenhuis bracket makes
$\fL_{\bcdot}=\Pi(C^{\bcdot+1}(\fn;\fn))$ into a Lie superalgebra that
can be interpreted as the centralizer of the exterior differential
considered as a vector field on the supermanifold $\hat M=(M,
\Omega(M))$ associated with the de Rham bundle on $M$.  A similar
bracket introduces structures of DG Lie superalgebra in
$\fL_{\bcdot}=\Pi(C^{\bcdot+1}(\fn;\fn))$ and
$\fl_{\bcdot}=\Pi(H^{\bcdot+1}(\fn;\fn))$ for any Lie superalgebra
$\fn$.  We explicitly describe (1) the algebras $\fl_0$ for the
maximal nilpotent subalgebra $\fn$ of any simple finite dimensional
Lie algebra $\fg$ and (2) the whole of $\fl_{\bcdot}$ for
$\fg=\fg_{2}$ which conjecturally has an archetypal structure, and in
the exceptional cases: $\fg=\fsl(2)$ and $\fsl(3)$; in $\fL_{\bcdot}$,
we also list the homologic elements (odd elements $x$ such that $[x,
x]=0$).  We observe that if the bracket in $\fn$ vanishes identically,
as is the case for Hermitean symmetric spaces $X=G/P$, the
differential $d$ is also zero and $\fl_{\bcdot}=\fvect(\dim 
\Pi(\fn))$, a simple Lie superalgebra.

We briefly review related results by Grozman, Penkov and
Serganova, Poletaeva, and Tolpygo.  We cite a powerful Premet's
theorem describing $H^{i}(\fn; N)$, where $N$ is an $\fn$-module which
is not a $\fg$-module.
\end{abstract}

\maketitle

\section*{\S 1. Introduction}

Though of 115 items of A.~L.~Onishchik's works listed today by
MathSciNet only a few are devoted to (co)homology and even these a few
deal with topological questions (as \cite{O}), rather than with
``linear algebra'' such as Lie algebra cohomology, the latter was one
of the main topics of the Vinberg-Onishchik seminar and Onishchik's
own and his students' studies some of which --- partly unpublished for
decades --- we review below.  It is at the above mentioned seminar
that one of us (DL) became intrigued by the Nijenhuis bracket and
deformation theory, a branch of which lately flourishes as the theory
of $L_{\infty}$-algebras.  The explicit definition of
$L_{\infty}$-algebras is irrelevant here and we refer the reader to
\cite{JS1}, \cite{JS2}, \cite{JS3}, \cite{HS}, \cite{M}.  Important
for us is that various DG (differential graded) Lie superalgebras are
examples of $L_{\infty}$-algebras.  There are, however, not many {\it
explicit} examples of $L_{\infty}$-algebras, even of DG Lie
superalgebras.

Here we give several explicit examples of DG Lie superalgebras.  The
naturalist trying to list the cohomology species analytically is set
back by the volume of the calculations.  That is why computer-aided
study, such as  \cite{Le}, the one we used here, or \cite{G2}-\cite{GL2}, is
indispensable.  We also list three tasks we intend to consider
shortly; our fourth task is to translate the theorems on
Lie algebras mentioned in what follows into the realm of Lie superalgebras.

\ssec{1.  Origins of the problem: Cohomology in physics and
mathematics} The electric charge and topological charges in gauge
theories are given by de Rham cohomology --- the pattern for
formulation of various (co)homology theories, in particular, for Lie
algebra cohomology, see \cite{GM}.  For further numerous (but still
covering only selected trends) examples, see, e.g., proceedings of the
International conference ``Cohomological Physics" \cite{CPh},
\cite{JS1}, \cite{AK}.  In short: cohomology is an important
invariant.  Here is one more application.

{\bf Curvature as cohomology}.  The case opposite in a sense to that
considered in \cite{B} is of particular interest.  Denote by $\fn_{i}$
the complement to a maximal parabolic subalgebra $\fp_{i}$ generated
by all simple coroots, both positive and negative, except the $i$-th
negative one, in any simple finite dimensional Lie algebra $\fg$.  If
$\fn_{i}$ is commutative, the elements of $H^2(\fn_{i}; \fg)$ can be
interpreted as the values of the generalization of the Riemann
curvature tensor at a point (see \cite{LPS}) and they ARE such values
if $\fg$ is the Lie algebra of $O(n)$ or $O(p,q)$.

More generally, let $P=P_{i_{1}\dots i_{k}}$ be any parabolic subgroup
(generated by all simple coroots, both positive and negative, except
the $i_{1}$-st, \dots and $i_{k}$-th negative ones) of a simple Lie
group $G$ and $N$ the complementary Lie group (i.e., $\fp\oplus
\fn=\fg$, where $\fp$, $\fn$ and $\fg$ are the Lie algebras of $P$,
$N$ and $G$, respectively).  Observe that if $\fn$ is not commutative,
the coset space $G/P$ is non-holonomic (here: is endowed with a
non-integrable distribution) and various cohomology of $\fn$ give
invariants of this non-holonomic manifold.  Some of such invariants
(the ones belonging to $H^2(\fn; \fg)$) were recently interpreted as
non-holonomic analogs of the curvature tensor, see \cite{L1}.  This
interpretation is actually a lucid expression of Wagner's description
of the nonholonomic analog of the curvature tensor performed in
components, see \cite{DG} and refs therein.

For Lie algebras $\fg$, the cohomology $H^2(\fn; \fg)$
was explicitly calculated until now only when $\fn$ is
commutative -- the cases of compact hermitian symmetric spaces
\cite{LPS}.  Now that we have Premet's theorem (below) we can use it
to extract an explicit description of $H^2(\fn; \fg)$ for
noncommutative $\fn$.  This is our nearest task.

The cohomology $H^i(\fn; \fg)$ for $i=1, 2$ constitute a part of the
Lie superalgebra $H^{\bcdot}(\fn; \fg)$ whose bracket is given by
point-wise bracketing of the skew-symmetric functions --- cocycles. 
To describe this Lie superalgebra is our second task.

The above mentioned bracket in $H^{\bcdot}(\fn; \fg)$ is given by a
scalar operator, and therefore is not so exciting as the bracket given
by a first order differential operator (on an appropriate
supermanifold) considered in this paper.

There are two major types of cohomology of $\fn$: with values in
$\fg$-modules and with other type of values.  Accordingly, there are
two types of results.

{\bf Example 1: The BWB theorem}.  The Borel-Weil-Bott
theorem (BWB)
states \cite{B} that, 
\vskip 0.2 cm

{\sl for any simple finite dimensional Lie algebra
$\fg$, its {\it maximal} nilpotent subalgebra $\fn_{\text{max}}$ and
any (finite dimensional) irreducible $\fg$-module $M$, the dimension
$H^i(\fn_{\text{max}}; M)$ does not depend on $M$ and is equal to the
cardinality of the set of elements of the Weyl group $W(\fg)$ of
length $i$.}
\vskip 0.2 cm

A generalization of the BWB theorem holds also for the nilpotent
subalgebras complementary to any parabolic subalgebra of $\fg$
(for references to ever more lucid and algebraic proofs of this algebraic
theorem, see \cite{T}).

When the Lie superalgebra theory started to develop, being boosted by
remarkable physical applications, to superize the BWB theorem was one
of the first problems.  It soon became clear that there is no hope to
get a neat super analog of the BWB theorem.  More precisely, for the
maximal nilpotent subalgebras, the theorems are relatively concise and
resemble the BWB theorem, though the notion of Weyl group becomes vague
and has several analogs in super setting (see \cite{Pe}, \cite{PS1},
\cite{PS2}).  For the ``opposite'' case of the complements to the
maximal parabolic superalgebras, the answer is rather complicated as
Poletaeva's results show \cite{Po1} (for the reader's convenience,
they are collected under one roof in \cite{Po2}).

Poletaeva's results, as well as \cite{LPS}, demonstrate inevitability
of computer-aided study (for more examples, see \cite{GL2}) of Lie
algebra (co)homology too complicated to deal with with bare hands;
especially so in super setting.

{\bf Example 2: From Kostant and Leger--Luks to Tolpygo and Premet}. 
There are plenty of $\fn$-modules which is not a $\fg$-modules, the
adjoint one is a most interesting.  In what follows we will briefly
write $C^{k}$ instead of $C^{k}(\fn; \fn)$ and similarly denote the
spaces of cocycles, coboundaries and cohomology ($Z^{k}$, $B^{k}$,
$H^{k}$, respectively).  We denote $H^{\bcdot}=\oplus H^{k}$, etc.

Kostant \cite{K} calculated $\dim H^1$; Leger and Luks \cite{LL}
calculated $\dim H^2$.  No general theorems on $H^{k}$ existed until
in mid-1970's Tolpygo \cite{T} computed $\dim H^{\bcdot}$ and
described the asymptotic behavior of $\dim H^k$ as $r \tto \infty$ for
$k<4$ and the {\it maximal} nilpotent subalgebra $\fn$ of classical
simple Lie algebra $\fg$ of series $\fsl(r)$, $\fo(r)$ and $\fsp(r)$.

To write supergravity equations, we need, at least, their left hand
sides, which are (parts of) nonholonomic curvature tensors with values
in $H^2(\fn; \fg)$ for some nilpotent subalgebras $\fn$ of
$\fsl(N|4)$, see \cite{GL1}, \cite{LPS}.  It is natural to compute
similar cohomology for nonholonomic manifolds as well, various flag
varieties, to start with.  DL discussed the matter with A.~Premet and
instead got the following beautiful general answer to a related
question \cite{Pr}.  The proof is regrettably still unpublished.

Hereafter the ground field is $\Cee$; let $\alpha _i$ denote the 
$i$th simple root of a simple Lie algebra $\fg$.

\begin{Theorem} {\em (Premet)} Let $\fg$ be a simple finite
dimensional Lie algebra, $\fp_{+}¥$ its parabolic subalgebra,
$\fp=\fp_{-}¥$ the opposite algebra and let $\fn$ be complementary to
$\fp_{+}¥$.  Let $E(\lambda) = L_{\lambda}$ be an irreducible (finite
dimensional) $\fg$-module with highest weight $\lambda$ such that $E
\simeq E^*$.  Let $V$ be a subspace in $E$ which is
$\fp_{-}$-invariant and contains $E_- = \mathop{\oplus}\limits_{\mu=
\sum k_i\alpha _i,\; k_i < 0}E_{\mu}$.  Let $P(V)$ be the set of
weights of $V$ and $w_{0}$ the longest element of the Weyl group and
let $\rho$ be the half sum of positive roots.  Then, for any $i < \rk~
\fg$, the following sequence is exact:
$$
\renewcommand{\arraystretch}{1.4}
\begin{array}{l}
    0 \tto \mathop{\oplus} \limits_{\{w\in W\mid
    l(w)=i-1, \, ww_0(\lambda)\not\in P(V)\}}
    E(-w(\lambda +\rho)+\rho)\tto H^{i-1}(\fn; E(\lambda)/V)\\
\tto  H^i(\fn; V)
\tto  \mathop{\oplus}
\limits_{\{w\in W\mid
    l(w)=i, \, ww_0(\lambda)\not\in P(V)\}}
E(-w(\lambda +\rho)+\rho)\tto  0.
\end{array}
$$
\end{Theorem}

It is very interesting to compute $\dim H^{\bcdot}$ and $\dim H^{k}$
for various subalgebras $\fn$ and all simple (finite dimensional, to
start with) Lie algebras $\fg$, and also Lie superalgebras.  But the
answer obtained in the form of dimensions of $H^{k}$, even if known,
is not quite satisfactory.

Indeed, the spaces $\fL_{\bcdot}=\Pi(C^{\bcdot+1})$ and
$\fl_{\bcdot}=\Pi(H^{\bcdot+1})$ are endowed with natural structures
of $\Zee$-graded Lie superalgebras (discovered by Nijenhuis in his,
with co-authors, studies of tensor invariants used in physics,
\cite{FN}, \cite{NR}).  Now observe that nobody presents graded
algebras by means of dimensions of their homogeneous components (e.g.,
we never introduce the polynomial algebra $S^{\bcdot}(\fg)$ as ``a
graded algebra whose degree $k$ component is of dimension
$\binom{n+k}{k}$''), we usually present algebras in terms of
generators and relations which is more graphic and more precise (the
homogeneous components of $U(\fg)$ have the same dimensions as those
of $S^{\bcdot}(\fg)$ for any deform of $\fg$).  Besides, given a Lie
(super)algebra, to ``{\it describe} it'' usually means to determine
its semisimple part and the radical.

The Lie superalgebras $\fL_{\bcdot}$ and $\fl_{\bcdot}$ are graded and
endowed with parity as follows: $\fL_j=\Pi(C^{j+1})$ and
$\fl_j=\Pi(H^{j+1})$.  In what follows we consider $\fl_{\bcdot}$. 
The part $\fl_{-1}$ is obvious (especially if $\fn$ is the {\it
maximal} nilpotent subalgebra of a simple $\fg$), and every $\fl_j$ is
an $\fl_0$-module.  In particular, if $\fn$ is a Lie algebra, then
$\fl_0$ is a Lie algebra.  The description of $\fl_0$ is a key step in
the description of $\fl_{\bcdot}$.  The first examples (calculation of
$\fl_{\bcdot}$ for $\fg=\fsl(2)$ and especially $\fsl(3)$) were very
encouraging.  Regrettably, for $\fg\neq\fsl(3)$ and its maximal
nilpotent subalgebra, the answer is not so neat.

\ssec{2.  Our result: examples of DG algebras and the list of their
homologic elements} The notion of $L_{\infty}$-algebras appeared as a
slackening of the notion of Lie algebra and a formulation of the
notions that vaguely lingered in various problems cf., e.g.,
\cite{JS2}, \cite{JS3}.  At the moment, there are not many explicit
examples of $L_{\infty}$-algebras, though it is known that Lie algebra
cohomology with values in the adjoint module, as well as the whole
space of cochains, and any DG algebra, are such examples.  Until this
paper, nothing was known (as far as we know) about the structure of
such DG algebras $\fl_{\bcdot}$.  To describe the whole of
$\fl_{\bcdot}$, in particular, for NON-maximal nilpotent sublagebras
$\fn$ of any simple (finite dimensional) Lie (super) algebra $\fg$ is
our third task.

Here we completely describe the Lie algebra $\fl_0=\Pi(H^1)$ for the
maximal nilpotent subalgebra of any simple (finite dimensional) Lie
algebra $\fg$ and also the whole of $\fl_{\bcdot}$ for $\fg=\fg_2$,
$\fsl(2)$ and $\fsl(3)$.  We conjecture that the structure of
$\fl_{\bcdot}$ is the same as for $\fg=\fg_2$ for any $\fg$, except
$\fsl(2)$ and $\fsl(3)$.  The formulation of the
result and conjecture (backed up with examples of $\fg=\fsl(n)$ for
$n\leq 5$) was molded with the help of computer-aided experiments
performed by one of us (\cite{Le}).

In numerous problems, the homologic elements in every Lie superalgebra
(odd elements $x$ such that $[x, x]=0$), or rather the homology they
determine, are vital, especially in deformation theory, see, e.g.,
\cite{Gs}, \cite{V2}, \cite{M}; they characterize algebroids \cite{V1}
and determine Lie algebras (\cite{L2}).  So we list them in
$\fl_{\bcdot}$ for $\fg=\fg_2$ and $\fsl(3)$.

\ssec{The Lie superalgebra structure on $\Pi(H^{\bcdot})$} Let $R^+$
be the set of positive roots, let the $e_\beta$, where $\beta\in R^+$,
be elements of the (weight) basis of $\fg$ (for example, the Chevalley
basis if $\fg$ is simple), let the $f^\beta$ be the elements of the
dual basis.  The basis of cochains is given by monomials $a, b, c\in
C^{\bcdot}$ of the form
$$
e_\alpha^i\otimes f_i^{\beta_{1}}\wedge ...\wedge
f_i^{\beta_{n_i}},~~~\text{where (by abuse of notation)}\;
i=a, b, c.
$$
Set
$$
a\cdot b=\sum_{k=1}^{n_{a}}
(-1)^{n_{a}-k}f_a^{\beta_{k}}(e_{\tilde\alpha}^b)e_{\alpha}^a\otimes
f_a^{\beta_{1}}\wedge ...\wedge \widehat{f_a^{\beta_{k}}}\wedge ... 
\wedge f_a^{\beta_{n_{a}}}\wedge f_b^{\tilde\beta_{n_{1}}}\wedge ... 
\wedge f_b^{\tilde\beta_{n_{b}}}
$$
and define the bracket on $\Pi(C^{\bcdot})$, where $\Pi$ is the shift
of parity functor, by the formula (hereafter $p(a)$ is the parity of
$a$ in $\Pi(C^{\bcdot})$, so, for Lie algebras $\fg$, it is equal to
$\deg a\pm 1$)
\begin{equation}
[a,b]=a\cdot b-(-1)^{p(a)p(b)}b\cdot a.
\end{equation}
It is subject to a direct verification that $(1)$ defines the Lie
superlagebra structure on $\fL_{\bcdot+1}=\Pi(C^{\bcdot+1})$, and
\begin{equation}
d[a,b]=[da,  b]+(-1)^{p(a)}[a, db].
\end{equation}
Therefore $\Pi(Z^{\bcdot+1})$ is a subalgebra of
$\Pi(C^{\bcdot+1})$ and $\Pi(B^{\bcdot+1})$ is an ideal in
$\Pi(Z^{\bcdot+1})$.  Hence, we have a DG Lie superalgebras structure
on $\fL_{\bcdot}:=\Pi(C^{\bcdot+1})$ and
$\fl_{\bcdot}:=\Pi(H^{\bcdot+1})$, cf., e.g., \cite{JS3}. The
differential in $\fl_{\bcdot}$ is zero and $\fl_{\bcdot}$ often is
very small, so $\Pi(Z^{\bcdot+1})$ might be more interesting than
$\fl_{\bcdot}$ in some questions, cf.  sec.  3.

For any maximal nilpotent subalgebra $\fn$ of any simple (finite
dimensional) Lie algebra $\fg$, the space $\fl_{-1}$ is, clearly,
one-dimensional.  

Let us describe the Lie algebra $\fl_{0}=\Pi(H^{1})$.  Select a basis
of $\fn$ consisting of root vectors $e_{\beta}$.  It is convenient to
select a basis in which the structure constants have the least
possible absolute values.  The Chevalley basis is (the only) such
basis and we select it for definiteness sake.  Let $n=\rk~\fg$.  The
$\Zee^n$-grading of $\fn$ by roots induces a $\Zee^n$-grading of
$C^{\bcdot}$ and we denote by $C^{i}_{\gamma}$ the subspace of $C^{i}$
of weight $\gamma$.

\ssec{2.1} To formulate our main result, we denote by $\mu$ the
(obviously, unique) maximal root in $\fg$, and by $l(\alpha)$
the level (i.e., the sum of coordinates with respect to the basis of
simple roots) of the root $\alpha$.  For $i=1, \ldots , n=\rk~\fg $, let
$$
m_i=\max\{m\in\Zee_{+}\mid \mu- k\alpha_i\in R^+ \text{ for any
}k=0, 1, \dots, m\}.
$$
Let $\alpha=\sum_{1\leq i\leq n} A_{i}(\alpha)\alpha_{i}$ be the
decomposition of a weight $\alpha$ with respect to the simple roots.

\begin{Theorem} If $\fg \neq \fsl(2)$, then $\dim H^1=2n$ {\em
(\cite{K})}, and for a basis one can take the cocycles
\begin{equation}
c_{i}=\sum_{\alpha\in R^+} A_{i}(\alpha)e_{\alpha}\otimes
f^{\alpha} 
\end{equation}
and
\begin{equation}
b_{i}=e_{\gamma_i}\otimes
f^{\alpha_i}, \quad\text{where }\; \gamma_i=\mu- m_i\alpha_i,
\end{equation}
where $i=1, \dots , n$. Further, for $\fg \neq \fsl(2), \fsl(3)$, we have
\begin{equation}
[c_i,c_j]=[b_i,b_j]=0,
\end{equation}
and
\begin{equation}
[c_i,b_j]=w_i(b_j)b_j,
\end{equation}
where $w(b)=(w_1(b),\dots,w_n(b))$ is the weight of the cochain $b$.
\end{Theorem}

\ssec{2.2. Proof}
We start with the following Lemma.

\ssbegin{2.2.1}{Lemma}{\em a)} $\dim H^1_{0}=n$ and for a basis one
can take the cocycles $(3)$.

{\em b)} $[c_i,c_j]=0$ for all $i$ and $j$.

\end{Lemma}

\begin{proof} a) Clearly, $B^1_{0}=0$, so $H^1_{0}=Z^1_{0}$.
Further, the cochain 
$$
c =\sum_{\alpha\in R^+} a(\alpha)e_{\alpha}\otimes
f^{\alpha}
$$ 
is a cocycle if and only if
\begin{equation}
a(\alpha+\beta)=a(\alpha)+a(\beta) \text{ whenever
$\alpha,\ \beta,\ \alpha+\beta\in R^+$.}
\end{equation}
It is evident that functions $A_i$ on $R^+$ satisfy $(6)$ and, by the
induction on the {\it level of the weight} (in our case, the sum of
the coefficients of the weight in the decomposition with respect to
simple roots, see \cite{FH}) we prove that the set $\{A_{i}\}_{i=1}^n$ forms a basis in the linear space of functions $a:
R^+\tto\Cee$ with property $(6)$.

b) is straightforward.
\end{proof}

\ssbegin{2.2.2}{Lemma}
The cochains $(4)$ are non-trivial linearly independent cocycles.
\end{Lemma}

\begin{proof}
First, observe that the $b_i$ are not co-boundaries because the
weight $\beta_i=\gamma_i-\alpha_i$ of $b_i$ is not a root, and hence
there are no co-boundaries of weight $\beta_i$ in $C^1$.

Let us verify that $db_i = 0$.  By definition, we have
$$
db_i=d(e_{\gamma_i}\otimes f^{\alpha_i})=de_{\gamma_i}\wedge f^{\alpha_i}-
e_{\mu}\otimes df^{\alpha_i}.
$$
It is clear that
$$
df^{\alpha_i}(e_{\alpha},e_{\beta})=f^{\alpha_i}([e_{\alpha},e_{\beta}])=
N_{\alpha\beta}f^{\alpha_i}(e_{\alpha+\beta})=0
$$
because $\alpha_i$ is a simple root and equality
$\alpha+\beta=\alpha_i$ is impossible for any positive roots $\alpha$
and $\beta$.  Further,
$$
de_{\gamma_i}(e_{\alpha})=[e_{\alpha},e_{\gamma_i}]=
N_{\alpha\gamma_i}e_{\alpha+\gamma_i}.
$$
The right hand side of this expression is nonzero only if
$\beta=\alpha+\gamma_i$ is a positive root.  Since $\mu$ is the
maximal root and $k\alpha_i$ is not a root for any simple root
$\alpha_i$ and any integer $k>1$, we see that
$$
de_{\gamma_i}=N_{\alpha_i\gamma_i}e_{\alpha_i+\gamma_i}\otimes f^{\alpha_i},
$$
where $N_{\alpha\beta}=0$ if $\alpha+\beta$ is not a root.
Hence,
$$
db_i=d(e_{\gamma_i}\otimes f^{\alpha_i})=de_{\gamma_i}\wedge
f^{\alpha_i}- e_{\gamma_i}\otimes df^{\alpha_i} =
N_{\alpha_i\gamma_i}e_{\gamma_i+\alpha_i}\otimes f^{\alpha_i}
\wedge f^{\alpha_i}=0.
$$
Since the weights of $b_i$ are distinct, they are linearly independent.
\end{proof}

Formula $(6)$ is a particular case of $(8)$, see sec. 2.4.3.

So we have $2n$ linearly independent nontrivial 1-cocycles; by 
Kostant's result (\cite{K}), they form a basis in $H^1$.  To complete
the proof of the theorem, we establish the following Lemma.

\ssbegin{2.2.3}{Lemma} If the level of $\mu$ is $\geq 5$, then
$[b_{i}, b_j]=0$ for all $i, j$.
\end{Lemma}
\begin{proof}
From the definition of the $b_i$ we see that
$$
[b_i,b_j]=f^{\alpha_i}(e_{\gamma_j})e_{\gamma_i}\otimes
f^{\alpha_j}- f^{\alpha_j}(e_{\gamma_i})e_{\gamma_j}\otimes
f^{\alpha_i}.
$$
As is known (\cite{FH}), in any simple finite dimensional Lie algebra
$\fg$, the ``arithmetical sequence'' of weights whose difference is a
positive root can not have more than four roots.  Hence, $m_i\leq 3$
and level of the root $\gamma_i=\mu-m_i\alpha_i$ is $\geq 2$.  So,
$\gamma_i$ is not a simple root and, in particular,
$\gamma_i\neq\alpha_j$.  Hence, $f^{\alpha_j}(e_{\gamma_i})=0$ for all
$i,j$.
\end{proof}

It is easy to see from the list of simple finite dimensional Lie
algebras that the level of $\mu$ is $<5$ only for the following $\fg$:
$$
A_1,\ A_2,\ A_3,\ A_4,\ B_2.
$$
Computer calculations show that (5) holds for all
algebras in the above list, except $A_1=\fsl(2)$ and $A_2=\fsl(3)$. \qed

\ssec{2.3} The ``exceptional" cases $\fsl(2)$ and $\fsl(3)$ are
described as follows.  \begin{Lemma} \underline{$\fsl(2)$:}
$\fl_{\bcdot}=\fl_{0}=\Span(c_1)$ is commutative;

\underline{$\fsl(3)$:} $\fl_{0}\simeq \fgl(2)$ with $[b_{1}, b_2]=c_2-c_1$.

Let $L^{n;\ c}$ be the irreducible $\fgl(2)$-module with the highest
weight $n$ with respect to $\fsl(2)$ and the value $c$ on the center,
$1_{2}\in \fgl(2)$. Then, as $\fl_{0}\simeq \fgl(2)$-modules,
$$
\fl_{-1}\simeq L^{0;2};\quad
\fl_{1}\simeq L^{2;-2}\oplus L^{1;-1} ;\quad
\fl_{2}\simeq L^{1;-3}.
$$
\end{Lemma}

\begin{proof} The statement about $\fsl(2)$ is obvious, that about
$\fsl(3)$ follows from the multiplication table 2.4.5. \end{proof}

\ssec{2.4. Lie superalgebra $\fl_{\bcdot}=\Pi(H^{\bcdot+1})$ for
$\fg = \fg_2$ and $\fsl(3)$}

\ssec{2.4.1. The case of $\fg=\fg_2$}
\begin{Theorem} For $\fg=\fg_2$, the basis of $\mathbf{H}^{\bcdot}$ is
given by the following list.  The multiplication table given in sec.
$2.4.2$ implies that $\fl_{\bcdot}$ is solvable and the
nonzero terms of the derived series  are
as follows (where $\fl_{(0)}=\fl$,
$\fl_{(i+1)}=[\fl,\fl_{(i)}]$~):
$$
\renewcommand{\arraystretch}{1.4}
\begin{array}{ll}
\fl_{(1)}=&\mathbf{H^0}\oplus \Span(h^1_3,h^1_4) \oplus
\mathbf{H^2}\oplus \mathbf{H^3}\oplus \mathbf{H^4}\oplus
\mathbf{H^5}\oplus\mathbf{H^6};\\
\fl_{(2)}=&\Span(h^1_3,\, h^1_4,\; h^2_2,\, h^2_4,\, h^2_6,\; h^3_2,\, h^3_3,
h^3_5,\, h^3_6,\, h^3_7,\, h^3_8,\\
&h^4_2,\, h^4_3,\, h^4_5,\, h^4_6,\, h^4_7,\;
h^5_2,\, h^5_3,\, h^5_4,\, h^5_5, \;h^6_1,\, h^6_2);\\
\fl_{(3)}= &\Span(h^2_4, h^6_2). \end{array}
$$
\end{Theorem}

\ssec{2.4.2. The basis of $\mathbf{H}^{\bcdot}$ for $\fg=\fg_{2}$}{$\ $}
$\mathbf{H}^{0}$:
$h^{0}_{1}=e_{ 2, 3}$

$\mathbf{H}^{1}$:

$h^{1}_{1}:=c_{1}=2e_{ 0, 1}\otimes f^{ 0, 1}-3e_{ 1,
0}\otimes f^{1, 0}-e_{ 1, 1}\otimes f^{ 1, 1}+e_{ 1, 2}\otimes f^{ 1,
2}+3e_{ 1, 3}\otimes f^{ 1, 3}$

$h^{1}_{2}:=c_{2}=e_{ 1, 0}\otimes f^{ 1, 0}+e_{ 1,
1}\otimes f^{ 1, 1}+e_{ 1, 2}\otimes f^{ 1, 2}+e_{ 1, 3}\otimes f^{ 1,
3}+2e_{ 2, 3}\otimes f^{ 2, 3}$

$h^{1}_{3}=e_{ 1, 3}\otimes f^{ 1, 0}$

$h^{1}_{4}=e_{ 2, 3}\otimes f^{ 0, 1}$

$\mathbf{H}^{2}$:

$h^{2}_{1}=-e_{ 0, 1}\otimes f^{ 1, 0}\wedge
f^{ 1, 1}+e_{ 1, 3}\otimes f^{ 1, 0}\wedge f^{ 2, 3}$

$h^{2}_{2}=-68e_{ 0, 1}\otimes f^{ 0, 1}\wedge f^{ 1, 1}+216e_{ 1,
0}\otimes f^{ 1, 0}\wedge f^{ 1, 1}+74e_{ 1, 1}\otimes f^{ 1, 0}\wedge
f^{ 1, 2}+39e_{ 1, 2}\otimes f^{ 1, 0}\wedge f^{ 1, 3}+31e_{ 1,
2}\otimes f^{ 1, 1}\wedge f^{ 1, 2}-9e_{ 1, 3}\otimes f^{ 0, 1}\wedge
f^{ 2, 3}+108e_{ 1, 3}\otimes f^{ 1, 1}\wedge f^{ 1, 3}+9e_{ 2,
3}\otimes f^{ 1, 1}\wedge f^{ 2, 3}$

$h^{2}_{3}=-e_{ 0, 1}\otimes f^{ 0, 1}\wedge
f^{ 1, 0}+e_{ 1, 1}\otimes f^{ 1, 0}\wedge f^{ 1, 1}+2e_{ 1, 2}\otimes
f^{ 1, 0}\wedge f^{ 1, 2}+3e_{ 1, 3}\otimes f^{ 1, 0}\wedge f^{ 1,
3}+3e_{ 2, 3}\otimes f^{ 1, 0}\wedge f^{ 2, 3}$

$h^{2}_{4}=e_{ 1, 3}\otimes f^{ 1, 0}\wedge f^{ 1, 1}$

$h^{2}_{5}=e_{ 1, 0}\otimes f^{ 0, 1}\wedge f^{ 1, 3}$

$h^{2}_{6}=11e_{ 1, 0}\otimes f^{ 0, 1}\wedge f^{ 1, 1}+14e_{ 1,
1}\otimes f^{ 0, 1}\wedge f^{ 1, 2}+9e_{ 1, 2}\otimes f^{ 0, 1}\wedge
f^{ 1, 3}+8e_{ 2, 3}\otimes f^{ 1, 2}\wedge f^{ 1, 3}$

$h^{2}_{7}=e_{ 1, 0}\otimes f^{ 0, 1}\wedge f^{ 1, 0}+e_{ 1, 1}\otimes
f^{ 0, 1}\wedge f^{ 1, 1}+e_{ 1, 2}\otimes f^{ 0, 1}\wedge f^{ 1,
2}+e_{ 1, 3}\otimes f^{ 0, 1}\wedge f^{ 1, 3}+2e_{ 2, 3}\otimes f^{ 0,
1}\wedge f^{ 2, 3}$

$\mathbf{H}^{3}$:

$h^{3}_{1}=-18e_{ 0, 1}\otimes f^{ 0, 1}\wedge f^{ 1, 1}\wedge f^{ 2,
3}-9e_{ 0, 1}\otimes f^{ 1, 0}\wedge f^{ 1, 2}\wedge f^{ 1, 3}+60e_{
1, 0}\otimes f^{ 1, 0}\wedge f^{ 1, 1}\wedge f^{ 2, 3}+21e_{ 1,
1}\otimes f^{ 1, 0}\wedge f^{ 1, 2}\wedge f^{ 2, 3}+10e_{ 1, 2}\otimes
f^{ 1, 0}\wedge f^{ 1, 3}\wedge f^{ 2, 3}+12e_{ 1, 2}\otimes f^{ 1,
1}\wedge f^{ 1, 2}\wedge f^{ 2, 3}+30e_{ 1, 3}\otimes f^{ 1, 1}\wedge
f^{ 1, 3}\wedge f^{ 2, 3}$

$h^{3}_{2}=13e_{ 0, 1}\otimes f^{ 1, 0}\wedge f^{ 1, 1}\wedge f^{ 1,
2}+6e_{ 1, 2}\otimes f^{ 1, 0}\wedge f^{ 1, 1}\wedge f^{ 2, 3}+9e_{ 1,
3}\otimes f^{ 1, 0}\wedge f^{ 1, 2}\wedge f^{ 2, 3}$

$h^{3}_{3}=e_{ 1, 3}\otimes f^{ 1, 0}\wedge f^{ 1, 1}\wedge f^{ 2, 3}$

$h^{3}_{4}=-217e_{ 0, 1}\otimes f^{ 0, 1}\wedge f^{ 1, 2}\wedge f^{ 1,
3}+124e_{ 1, 0}\otimes f^{ 0, 1}\wedge f^{ 1, 1}\wedge f^{ 2, 3}+62e_{
1, 0}\otimes f^{ 1, 0}\wedge f^{ 1, 2}\wedge f^{ 1, 3}+183e_{ 1,
1}\otimes f^{ 0, 1}\wedge f^{ 1, 2}\wedge f^{ 2, 3}+28e_{ 1, 1}\otimes
f^{ 1, 1}\wedge f^{ 1, 2}\wedge f^{ 1, 3}+126e_{ 1, 2}\otimes f^{ 0,
1}\wedge f^{ 1, 3}\wedge f^{ 2, 3}+127e_{ 2, 3}\otimes f^{ 1, 2}\wedge
f^{ 1, 3}\wedge f^{ 2, 3}$

$h^{3}_{5}=-3e_{ 0, 1}\otimes f^{ 0, 1}\wedge f^{ 1, 0}\wedge f^{ 1,
1}-3e_{ 1, 2}\otimes f^{ 1, 0}\wedge f^{ 1, 1}\wedge f^{ 1, 2}+e_{ 1,
3}\otimes f^{ 0, 1}\wedge f^{ 1, 0}\wedge f^{ 2, 3}-5e_{ 1, 3}\otimes
f^{ 1, 0}\wedge f^{ 1, 1}\wedge f^{ 1, 3}-4e_{ 2, 3}\otimes f^{ 1,
0}\wedge f^{ 1, 1}\wedge f^{ 2, 3}$

$h^{3}_{6}=e_{ 1, 0}\otimes f^{ 0, 1}\wedge f^{1, 2}\wedge f^{1, 3}$

$h^{3}_{7}=e_{ 1, 0}\otimes f^{ 0, 1}\wedge f^{ 1, 1}\wedge f^{ 1,
3}+2e_{ 1, 1}\otimes f^{ 0, 1}\wedge f^{1, 2}\wedge f^{1, 3}$

$h^{3}_{8}=28e_{ 1, 0}\otimes f^{ 0, 1}\wedge f^{ 1, 0}\wedge f^{ 1,
3}-3e_{ 1, 0}\otimes f^{ 0, 1}\wedge f^{ 1, 1}\wedge f^{ 1, 2}+19e_{
1, 1}\otimes f^{ 0, 1}\wedge f^{ 1, 1}\wedge f^{ 1, 3}+19e_{ 1,
2}\otimes f^{ 0, 1}\wedge f^{ 1, 2}\wedge f^{ 1, 3}-37e_{ 2, 3}\otimes
f^{ 0, 1}\wedge f^{ 1, 3}\wedge f^{ 2, 3}$

$\mathbf{H}^{4}$:

$h^{4}_{1}=e_{ 0, 1}\otimes f^{ 1, 0}\wedge f^{ 1, 1}\wedge f^{ 1,
2}\wedge f^{ 2, 3}$

$h^{4}_{2}=-4e_{ 0, 1}\otimes f^{ 0, 1}\wedge f^{ 1, 0}\wedge f^{ 1,
2}\wedge f^{ 2, 3}+35e_{ 0, 1}\otimes f^{ 1, 0}\wedge f^{ 1, 1}\wedge
f^{ 1, 2}\wedge f^{ 1, 3}-31e_{ 1, 1}\otimes f^{ 1, 0}\wedge f^{ 1,
1}\wedge f^{ 1, 2}\wedge f^{ 2, 3}-18e_{ 1, 2}\otimes f^{ 1, 0}\wedge
f^{ 1, 1}\wedge f^{ 1, 3}\wedge f^{ 2, 3}-27e_{ 1, 3}\otimes f^{ 1,
0}\wedge f^{ 1, 2}\wedge f^{ 1, 3}\wedge f^{ 2, 3}$

$h^{4}_{3}=e_{ 0, 1}\otimes f^{ 0, 1}\wedge f^{ 1, 0}\wedge f^{ 1,
1}\wedge f^{ 2, 3}+e_{ 1, 2}\otimes f^{ 1, 0}\wedge f^{ 1, 1}\wedge
f^{ 1, 2}\wedge f^{ 2, 3}+2e_{ 1, 3}\otimes f^{ 1, 0}\wedge f^{ 1,
1}\wedge f^{ 1, 3}\wedge f^{ 2, 3}$

$h^{4}_{4}=e_{ 1, 0}\otimes f^{ 0, 1}\wedge f^{ 1, 2}\wedge f^{ 1,3}
\wedge f^{ 2, 3}$

$h^{4}_{5}=e_{ 1, 0}\otimes f^{ 0, 1}\wedge f^{ 1, 1}\wedge f^{ 1,
3}\wedge f^{ 2, 3}+2e_{ 1, 1}\otimes f^{ 0, 1}\wedge f^{ 1, 2}\wedge
f^{ 1, 3}\wedge f^{ 2, 3}$

$h^{4}_{6}=11e_{ 0, 1}\otimes f^{ 0, 1}\wedge f^{ 1, 1}\wedge f^{ 1,
2}\wedge f^{ 1, 3}-4e_{ 1, 0}\otimes f^{ 0, 1}\wedge f^{ 1, 0}\wedge
f^{ 1, 2}\wedge f^{ 2, 3}-4e_{ 1, 0}\otimes f^{ 1, 0}\wedge f^{ 1,
1}\wedge f^{ 1, 2}\wedge f^{ 1, 3}-7e_{ 1, 1}\otimes f^{ 0, 1}\wedge
f^{ 1, 1}\wedge f^{ 1, 2}\wedge f^{ 2, 3}-9e_{ 1, 2}\otimes f^{ 0,
1}\wedge f^{ 1, 1}\wedge f^{ 1, 3}\wedge f^{ 2, 3}-5e_{ 1, 3}\otimes
f^{ 0, 1}\wedge f^{ 1, 2}\wedge f^{ 1, 3}\wedge f^{ 2, 3}+e_{ 2,
3}\otimes f^{ 1, 1}\wedge f^{ 1, 2}\wedge f^{ 1, 3}\wedge f^{ 2, 3}$

$h^{4}_{7}=e_{ 1, 0}\otimes f^{ 0, 1}\wedge f^{ 1, 0}\wedge f^{ 1,
2}\wedge f^{ 1, 3}+e_{ 1, 1}\otimes f^{ 0, 1}\wedge f^{ 1, 1}\wedge
f^{ 1, 2}\wedge f^{ 1, 3}+e_{ 2, 3}\otimes f^{ 0, 1}\wedge f^{ 1,
2}\wedge f^{ 1, 3}\wedge f^{ 2, 3}$

$\mathbf{H}^{5}$:

$h^{5}_{1}=e_{ 0, 1}\otimes f^{ 1, 0}\wedge f^{ 1, 1}\wedge f^{ 1,
2}\wedge f^{ 1, 3}\wedge f^{ 2, 3}$

$h^{5}_{2}=-e_{ 0, 1}\otimes f^{ 0, 1}\wedge f^{ 1, 1}\wedge f^{ 1,
2}\wedge f^{ 1, 3}\wedge f^{ 2, 3}+e_{ 1, 0}\otimes f^{ 1, 0}\wedge
f^{ 1, 1}\wedge f^{ 1, 2}\wedge f^{ 1, 3}\wedge f^{ 2, 3}$

$h^{5}_{3}=e_{ 0, 1}\otimes f^{ 0, 1}\wedge f^{ 1, 0}\wedge f^{ 1,
1}\wedge f^{ 1, 2}\wedge f^{ 2, 3}-e_{ 1, 3}\otimes f^{ 1, 0}\wedge
f^{ 1, 1}\wedge f^{ 1, 2}\wedge f^{ 1, 3}\wedge f^{ 2, 3}$

$h^{5}_{4}=e_{ 1, 0}\otimes f^{ 0, 1}\wedge f^{ 1, 1}\wedge f^{ 1,
2}\wedge f^{ 1, 3}\wedge f^{ 2, 3}$

$h^{5}_{5}=e_{ 1, 0}\otimes f^{ 0, 1}\wedge f^{ 1, 0}\wedge f^{ 1,
2}\wedge f^{ 1, 3}\wedge f^{ 2, 3}+e_{ 1, 1}\otimes f^{ 0, 1}\wedge
f^{ 1, 1}\wedge f^{ 1, 2}\wedge f^{ 1, 3}\wedge f^{ 2, 3}$

$\mathbf{H}^{6}$:

$h^{6}_{1}=e_{ 0, 1}\otimes f^{ 0, 1}\wedge f^{ 1, 0}\wedge f^{ 1,
1}\wedge f^{ 1, 2}\wedge f^{ 1, 3}\wedge f^{ 2, 3}$

$h^{6}_{2}=e_{ 1, 0}\otimes f^{ 0, 1}\wedge f^{ 1, 0}\wedge f^{ 1,
1}\wedge f^{ 1, 2}\wedge f^{ 1, 3}\wedge f^{ 2, 3}$

\ssec{2.4.3.  The multiplication table in $\mathbf{H}^{\bcdot}$}
Observe that the Cartan subalgebra $\fh$ of $\fg$ naturally acts on
$\fn$ and $C^{\bcdot}$.  The image of $\fh$ in $C^{1}=\fgl(\fn)$ is
precisely $\Span(c_{i}\mid i=1, \dots , \rk\;\fg)$, so we will save
several lines in the multiplication tables $H^1\times H^k$ here and in
sec.  2.4.5 by replacing them with a generalization of formula (6):
\begin{equation}
[c_{i}, c]=w_{i}(c)c\text{ for any chain $c$ of weight $w(c)$.}
\end{equation}

$\mathbf{H}^{0}\times\mathbf{H}^{1}$:

\begin{center}
\renewcommand{\arraystretch}{1.4}
\begin{tabular}{|c|c|c|c|c|}
\hline
   & $h^{1}_{1}$& $h^{1}_{2}$& $h^{1}_{3}$& $h^{1}_{4}$\\
\hline
$h^{0}_{1}$& $0$& $-2 h^{0}_{1}$& $0$& $0$\\
\hline
\end{tabular}\end{center}

$\mathbf{H}^{0}\times\mathbf{H}^{2}$:

\begin{center}
\renewcommand{\arraystretch}{1.4}
\begin{tabular}{|c|c|c|c|c|c|c|c|}
\hline
   & $h^{2}_{1}$& $h^{2}_{2}$& $h^{2}_{3}$& $h^{2}_{4}$& $h^{2}_{5}$&
   $h^{2}_{6}$& $h^{2}_{7}$\\
\hline
$h^{0}_{1}$& $h^{1}_{3}$& $0$& $0$& $0$& $0$& $0$& $2 h^{1}_{4}$\\
\hline
\end{tabular}\end{center}

$\mathbf{H}^{0}\times\mathbf{H}^{3}$:

\begin{center}
\renewcommand{\arraystretch}{1.4}
\begin{tabular}{|c|c|c|c|c|c|c|c|c|}
\hline
   & $h^{3}_{1}$& $h^{3}_{2}$& $h^{3}_{3}$& $h^{3}_{4}$& $h^{3}_{5}$&
   $h^{3}_{6}$& $h^{3}_{7}$& $h^{3}_{8}$\\
\hline
$h^{0}_{1}$& $-\frac{5}{18} h^{2}_{2}$& $0$& $-h^{2}_{4}$&
$-\frac{434}{33} h^{2}_{6}$& $0$& $0$& $0$& $0$\\
\hline
\end{tabular}\end{center}

$\mathbf{H}^{0}\times\mathbf{H}^{4}$:

\begin{center}
    \renewcommand{\arraystretch}{1.4}
\begin{tabular}{|c|c|c|c|c|c|c|c|}
\hline
   & $h^{4}_{1}$& $h^{4}_{2}$& $h^{4}_{3}$& $h^{4}_{4}$& $h^{4}_{5}$&
   $h^{4}_{6}$& $h^{4}_{7}$\\
\hline
$h^{0}_{1}$& $\frac{1}{22} h^{3}_{2}$& $0$& $-\frac{4}{15} h^{3}_{5}$&
$h^{3}_{6}$& $h^{3}_{7}$& $0$& $0$\\
\hline
\end{tabular}\end{center}

$\mathbf{H}^{0}\times\mathbf{H}^{5}$:

\begin{center}
\renewcommand{\arraystretch}{1.4}
\begin{tabular}{|c|c|c|c|c|c|}
\hline
   & $h^{5}_{1}$& $h^{5}_{2}$& $h^{5}_{3}$& $h^{5}_{4}$& $h^{5}_{5}$\\
\hline
$h^{0}_{1}$& $-\frac{1}{93} h^{4}_{2}$& $\frac{5}{103} h^{4}_{6}$& $0$&
$0$& $-\frac{2}{3} h^{4}_{7}$\\
\hline
\end{tabular}\end{center}

$\mathbf{H}^{0}\times\mathbf{H}^{6}=0$:

$\mathbf{H}^{1}\times\mathbf{H}^{1}=0$ (except for (8))

$\mathbf{H}^{1}\times\mathbf{H}^{2}$:

\begin{center}\renewcommand{\arraystretch}{1.4}
\begin{tabular}{|c|c|c|c|c|c|c|c|}
\hline
   & $h^{2}_{1}$& $h^{2}_{2}$& $h^{2}_{3}$& $h^{2}_{4}$& $h^{2}_{5}$&
   $h^{2}_{6}$& $h^{2}_{7}$\\
\hline
$h^{1}_{3}$& $0$& $324 h^{2}_{4}$& $0$& $0$& $0$& $0$& $0$\\
\hline
$h^{1}_{4}$& $0$& $0$& $0$& $0$& $0$& $0$& $0$\\
\hline
\end{tabular}\end{center}

$\mathbf{H}^{1}\times\mathbf{H}^{3}$:

\begin{center}\renewcommand{\arraystretch}{1.4}
\begin{tabular}{|c|c|c|c|c|c|c|c|c|}
\hline
   & $h^{3}_{1}$& $h^{3}_{2}$& $h^{3}_{3}$& $h^{3}_{4}$& $h^{3}_{5}$&
   $h^{3}_{6}$& $h^{3}_{7}$& $h^{3}_{8}$\\
\hline
$h^{1}_{3}$& $90 h^{3}_{3}$& $0$& $0$& $0$& $0$& $0$& $0$& $0$\\
\hline
$h^{1}_{4}$& $0$& $0$& $0$& $0$& $0$& $0$& $0$& $0$\\
\hline
\end{tabular}\end{center}

$\mathbf{H}^{1}\times\mathbf{H}^{4}$:

\begin{center}\renewcommand{\arraystretch}{1.4}
\begin{tabular}{|c|c|c|c|c|c|c|c|}
\hline
   & $h^{4}_{1}$& $h^{4}_{2}$& $h^{4}_{3}$& $h^{4}_{4}$& $h^{4}_{5}$&
   $h^{4}_{6}$& $h^{4}_{7}$\\
\hline
$h^{1}_{3}$& $0$& $0$& $0$& $-\frac{3}{103} h^{4}_{6}$& $0$& $0$& $0$\\
\hline
$h^{1}_{4}$& $0$& $0$& $0$& $0$& $0$& $0$& $0$\\
\hline
\end{tabular}\end{center}

$\mathbf{H}^{1}\times\mathbf{H}^{5}=0$ (except for (8))

$\mathbf{H}^{1}\times\mathbf{H}^{6}=0$ (except for (8))

$\mathbf{H}^{2}\times\mathbf{H}^{2}$:

\begin{center}\renewcommand{\arraystretch}{1.4}
\begin{tabular}{|c|c|c|c|c|c|c|c|}
\hline
   & $h^{2}_{1}$& $h^{2}_{2}$& $h^{2}_{3}$& $h^{2}_{4}$& $h^{2}_{5}$&
   $h^{2}_{6}$& $h^{2}_{7}$\\
\hline
$h^{2}_{1}$&0 & $-324 h^{3}_{3}$& $0$& $0$& $0$& $0$& $-\frac{2}{5}
h^{3}_{5}$\\
\hline
$h^{2}_{2}$&- &0 & $144 h^{3}_{5}$& $0$& $108 h^{3}_{7}$& $0$& $0$\\
\hline
$h^{2}_{3}$&- &- &0 & $0$& $\frac{10}{103} h^{3}_{8}$& $0$& $0$\\
\hline
$h^{2}_{4}$&- &- &- &0 & $0$& $0$& $0$\\
\hline
$h^{2}_{5}$&- &- &- &- &0 & $0$& $0$\\
\hline
$h^{2}_{6}$&- &- &- &- &- &0 & $0$\\
\hline
$h^{2}_{7}$&- &- &- &- &- &- &0 \\
\hline
\end{tabular}\end{center}

$\mathbf{H}^{2}\times\mathbf{H}^{3}$:

\begin{center}\renewcommand{\arraystretch}{1.4}
\begin{tabular}{|c|c|c|c|c|c|c|c|c|}
\hline
   & $h^{3}_{1}$& $h^{3}_{2}$& $h^{3}_{3}$& $h^{3}_{4}$& $h^{3}_{5}$&
   $h^{3}_{6}$& $h^{3}_{7}$& $h^{3}_{8}$\\
\hline
$h^{2}_{1}$& $0$& $0$& $0$& $-7 h^{4}_{2}$& $0$& $-\frac{8}{103}
h^{4}_{6}$& $0$& $0$\\
\hline
$h^{2}_{2}$& $\frac{4320}{31} h^{4}_{2}$& $0$& $0$&
$-\frac{140616}{103} h^{4}_{6}$& $0$& $0$& $0$& $0$\\
\hline
$h^{2}_{3}$& $-150 h^{4}_{3}$& $0$& $0$& $0$& $0$& $\frac{10}{3}
h^{4}_{7}$& $0$& $0$\\
\hline
$h^{2}_{4}$& $0$& $0$& $0$& $0$& $0$& $0$& $0$& $0$\\
\hline
$h^{2}_{5}$& $30 h^{4}_{5}$& $-\frac{110}{103} h^{4}_{6}$& $0$& $0$&
$0$& $0$& $0$& $0$\\
\hline
$h^{2}_{6}$& $\frac{1980}{103} h^{4}_{6}$& $0$& $0$& $0$& $0$& $0$&
$0$& $0$\\
\hline
$h^{2}_{7}$& $0$& $0$& $0$& $-434 h^{4}_{7}$& $0$& $0$& $0$& $0$\\
\hline
\end{tabular}\end{center}

$\mathbf{H}^{2}\times\mathbf{H}^{4}$:

\begin{center}\renewcommand{\arraystretch}{1.4}
\begin{tabular}{|c|c|c|c|c|c|c|c|}
\hline
   & $h^{4}_{1}$& $h^{4}_{2}$& $h^{4}_{3}$& $h^{4}_{4}$& $h^{4}_{5}$&
   $h^{4}_{6}$& $h^{4}_{7}$\\
\hline
$h^{2}_{1}$& $0$& $0$& $0$& $h^{5}_{2}$& $0$& $0$& $0$\\
\hline
$h^{2}_{2}$& $0$& $0$& $0$& $432 h^{5}_{4}$& $0$& $0$& $0$\\
\hline
$h^{2}_{3}$& $0$& $0$& $0$& $5 h^{5}_{5}$& $0$& $0$& $0$\\
\hline
$h^{2}_{4}$& $0$& $0$& $0$& $0$& $0$& $0$& $0$\\
\hline
$h^{2}_{5}$& $h^{5}_{2}$& $0$& $0$& $0$& $0$& $0$& $0$\\
\hline
$h^{2}_{6}$& $0$& $0$& $0$& $0$& $0$& $0$& $0$\\
\hline
$h^{2}_{7}$& $-3 h^{5}_{3}$& $0$& $0$& $0$& $0$& $0$& $0$\\
\hline
\end{tabular}\end{center}

$\mathbf{H}^{2}\times\mathbf{H}^{5}$:

\begin{center}\renewcommand{\arraystretch}{1.4}
\begin{tabular}{|c|c|c|c|c|c|}
\hline
   & $h^{5}_{1}$& $h^{5}_{2}$& $h^{5}_{3}$& $h^{5}_{4}$& $h^{5}_{5}$\\
\hline
$h^{2}_{1}$& $0$& $0$& $0$& $0$& $2 h^{6}_{1}$\\
\hline
$h^{2}_{2}$& $0$& $0$& $0$& $0$& $-432 h^{6}_{2}$\\
\hline
$h^{2}_{3}$& $0$& $-10 h^{6}_{1}$& $0$& $10 h^{6}_{2}$& $0$\\
\hline
$h^{2}_{4}$& $0$& $0$& $0$& $0$& $0$\\
\hline
$h^{2}_{5}$& $0$& $0$& $-2 h^{6}_{2}$& $0$& $0$\\
\hline
$h^{2}_{6}$& $0$& $0$& $0$& $0$& $0$\\
\hline
$h^{2}_{7}$& $-6 h^{6}_{1}$& $-6 h^{6}_{2}$& $0$& $0$& $0$\\
\hline
\end{tabular}\end{center}

$\mathbf{H}^{3}\times\mathbf{H}^{3}$:

\begin{center}\renewcommand{\arraystretch}{1.4}
\begin{tabular}{|c|c|c|c|c|c|c|c|c|}
\hline
   & $h^{3}_{1}$& $h^{3}_{2}$& $h^{3}_{3}$& $h^{3}_{4}$& $h^{3}_{5}$&
   $h^{3}_{6}$& $h^{3}_{7}$& $h^{3}_{8}$\\
\hline
$h^{3}_{1}$&- & $0$& $0$& $13020 h^{5}_{2}$& $0$& $120 h^{5}_{4}$&
$0$& $0$\\
\hline
$h^{3}_{2}$&- &- & $0$& $0$& $0$& $0$& $0$& $0$\\
\hline
$h^{3}_{3}$&- &- &- & $0$& $0$& $0$& $0$& $0$\\
\hline
$h^{3}_{4}$&- &- &- &- & $0$& $0$& $0$& $0$\\
\hline
$h^{3}_{5}$&- &- &- &- &- & $0$& $0$& $0$\\
\hline
$h^{3}_{6}$&- &- &- &- &- &- & $0$& $0$\\
\hline
$h^{3}_{7}$&- &- &- &- &- &- &- & $0$\\
\hline
$h^{3}_{8}$&- &- &- &- &- &- &- &- \\
\hline
\end{tabular}\end{center}

$\mathbf{H}^{3}\times\mathbf{H}^{4}$:

\begin{center}\renewcommand{\arraystretch}{1.4}
\begin{tabular}{|c|c|c|c|c|c|c|c|}
\hline
   & $h^{4}_{1}$& $h^{4}_{2}$& $h^{4}_{3}$& $h^{4}_{4}$& $h^{4}_{5}$&
   $h^{4}_{6}$& $h^{4}_{7}$\\
\hline
$h^{3}_{1}$& $0$& $0$& $0$& $0$& $0$& $0$& $180 h^{6}_{2}$\\
\hline
$h^{3}_{2}$& $0$& $0$& $0$& $0$& $0$& $0$& $0$\\
\hline
$h^{3}_{3}$& $0$& $0$& $0$& $0$& $0$& $0$& $0$\\
\hline
$h^{3}_{4}$& $0$& $0$& $-868 h^{6}_{1}$& $0$& $0$& $0$& $0$\\
\hline
$h^{3}_{5}$& $0$& $0$& $0$& $15 h^{6}_{2}$& $0$& $0$& $0$\\
\hline
$h^{3}_{6}$& $0$& $0$& $4 h^{6}_{2}$& $0$& $0$& $0$& $0$\\
\hline
$h^{3}_{7}$& $0$& $0$& $0$& $0$& $0$& $0$& $0$\\
\hline
$h^{3}_{8}$& $-103 h^{6}_{1}$& $0$& $0$& $0$& $0$& $0$& $0$\\
\hline
\end{tabular}\end{center}

For the (boring) details of the proof and the (lively) program we used, see
\cite{Le}.

\ssec{2.4.4. The basis of $\mathbf{H}^{\bcdot}$ for $\fg=\fsl(3)$}{$\ $}
$\mathbf{H}^{0}:$
$h^{0}_{1}=e_{1, 1}$

$\mathbf{H}^{1}$ (where $c_{i}=E_{ii}$, $h^{1}_{4}$ and
$h^{1}_{1}$ are the raising and lowing operators of $\fsl(2)$):

$h^{1}_{1}=e_{ 0, 1}\otimes f^{ 1, 0}$

$h^{1}_{2}:=c_{1}-c_{2}=-e_{ 0, 1}\otimes f^{ 0, 1}+e_{ 1, 0}\otimes f^{ 1, 0}$

$h^{1}_{3}:=c_{1}+c_{2}=e_{ 0, 1}\otimes f^{ 0, 1}+e_{ 1, 0}\otimes f^{ 1, 0}+2e_{
1, 1}\otimes f^{ 1, 1}$

$h^{1}_{4}=e_{ 1, 0}\otimes f^{ 0, 1}$

$\mathbf{H}^{2}:$

$h^{2}_{1}=e_{ 0, 1}\otimes f^{ 1, 0}\wedge f^{ 1, 1}$

$h^{2}_{2}=-e_{ 0, 1}\otimes f^{ 0, 1}\wedge f^{ 1, 1}+e_{ 1,
0}\otimes f^{ 1, 0}\wedge f^{ 1, 1}$

$h^{2}_{3}=e_{ 0, 1}\otimes f^{ 0, 1}\wedge f^{ 1, 0}-e_{ 1, 1}\otimes
f^{ 1, 0}\wedge f^{ 1, 1}$

$h^{2}_{4}=e_{ 1, 0}\otimes f^{ 0, 1}\wedge f^{ 1, 1}$

$h^{2}_{5}=e_{ 1, 0}\otimes f^{ 0, 1}\wedge f^{ 1, 0}+e_{ 1, 1}\otimes
f^{ 0, 1}\wedge f^{ 1, 1}$

$\mathbf{H}^{3}:$

$h^{3}_{1}=e_{ 0, 1}\otimes f^{ 0, 1}\wedge f^{ 1, 0}\wedge f^{ 1, 1}$

$h^{3}_{2}=e_{ 1, 0}\otimes f^{ 0, 1}\wedge f^{ 1, 0}\wedge f^{ 1, 1}$

\ssec{2.4.5. The multiplication table in $\fl_{\bcdot}$ for $\fg=\fsl(3)$}
{}~{}

$\mathbf{H}^{0}\times\mathbf{H}^{1}$:

\begin{center}
\begin{tabular}{|c|c|c|c|c|}
\hline
   & $h^{1}_{1}$& $h^{1}_{2}$& $h^{1}_{3}$& $h^{1}_{4}$\\
\hline
$h^{0}_{1}$& $0$& $0$& $-2 h^{0}_{1}$& $0$\\
\hline
\end{tabular}\end{center}

$\mathbf{H}^{0}\times\mathbf{H}^{2}$:

\begin{center}
\begin{tabular}{|c|c|c|c|c|c|}
\hline
   & $h^{2}_{1}$& $h^{2}_{2}$& $h^{2}_{3}$& $h^{2}_{4}$& $h^{2}_{5}$\\
\hline
$h^{0}_{1}$& $h^{1}_{1}$& $h^{1}_{2}$& $0$& $h^{1}_{4}$& $0$\\
\hline
\end{tabular}\end{center}

$\mathbf{H}^{0}\times\mathbf{H}^{3}$:

\begin{center}
\begin{tabular}{|c|c|c|}
\hline
   & $h^{3}_{1}$& $h^{3}_{2}$\\
\hline
$h^{0}_{1}$& $-\frac12 h^{2}_{3}$& $-\frac12 h^{2}_{5}$\\
\hline
\end{tabular}\end{center}

$\mathbf{H}^{1}\times\mathbf{H}^{1}$:

\begin{center}
\begin{tabular}{|c|c|c|}
\hline
   & $h^{1}_{1}$&$h^{1}_{4}$\\
\hline
$h^{1}_{1}$& $0$& $-h^{1}_{2}$\\
\hline
$h^{1}_{4}$& $h^{1}_{2}$&$0$\\
\hline
\end{tabular}\end{center}

$\mathbf{H}^{1}\times\mathbf{H}^{2}$:

\begin{center}
\begin{tabular}{|c|c|c|c|c|c|}
\hline
   & $h^{2}_{1}$& $h^{2}_{2}$& $h^{2}_{3}$& $h^{2}_{4}$& $h^{2}_{5}$\\
\hline
$h^{1}_{1}$& $0$& $2 h^{2}_{1}$& $0$& $-h^{2}_{2}$& $h^{2}_{3}$\\
\hline
$h^{1}_{4}$& $h^{2}_{2}$& $-2 h^{2}_{4}$& $h^{2}_{5}$& $0$& $0$\\
\hline
\end{tabular}\end{center}

$\mathbf{H}^{1}\times\mathbf{H}^{3}$:

\begin{center}
\begin{tabular}{|c|c|c|}
\hline
   & $h^{3}_{1}$& $h^{3}_{2}$\\
\hline
$h^{1}_{1}$& $0$& $h^{3}_{1}$\\
\hline
$h^{1}_{4}$& $h^{3}_{2}$& $0$\\
\hline
\end{tabular}\end{center}

$\mathbf{H}^{2}\times\mathbf{H}^{2}$:

\begin{center}
\begin{tabular}{|c|c|c|c|c|c|}
\hline
   & $h^{2}_{1}$& $h^{2}_{2}$& $h^{2}_{3}$& $h^{2}_{4}$& $h^{2}_{5}$\\
\hline
$h^{2}_{1}$& $0$& $0$& $0$& $0$& $-2 h^{3}_{1}$\\
\hline
$h^{2}_{2}$& $0$& $0$& $2 h^{3}_{1}$& $0$& $-2 h^{3}_{2}$\\
\hline
$h^{2}_{3}$& $0$& $2 h^{3}_{1}$& $0$& $-2 h^{3}_{2}$& $0$\\
\hline
$h^{2}_{4}$& $0$& $0$& $-2 h^{3}_{2}$& $0$& $0$\\
\hline
$h^{2}_{5}$& $-2 h^{3}_{1}$& $-2 h^{3}_{2}$& $0$& $0$& $0$\\
\hline
\end{tabular}\end{center}

\ssec{2.5.1. The homologic elements in $\fl_{\bcdot}$ for $\fg=\fg_2$}
Let $x\in\fl_{\bcdot}$ be a homologic element (odd and
such that $[x, x]=0$).  Then $x$ is of the
form
\begin{equation}
x=x^0+x^2+x^4+x^6,\; \text{ where }\; x^i=\sum_j k^i_j h^i_j
\in\mathbf{H}^i.
\end{equation}

a) $x^0=0$.  In this case the condition
$[x,x]=0$ takes the form
$$
[x^2,x^2]=0\quad \text{ and }\quad [x^2,x^4]=0.
$$
Let us find first the form of $x^2$.  From the multiplication table we
deduce that $[x^2,x^2]=0$ if and only if the support of (the indices
of) nonzero coefficients in the sum $x^2=\sum_j k^2_j h^2_j$ --- a
subset of the set $\{1, 2, \dots, 7\}$ --- is of the form $A\cup B$,
where $A\subset \{4, 6\}$ and $B$ is only one of following
$$
\emptyset; \; \{1\}; \;\{2\}; \;\{3\}; \;\{5\}; \;\{7\}; \;\{1,3\};
\;\{1,5\}; \;\{2,7\}; \;\{3,7\} ; \;\{5,7\}.
$$
We have $[x^2,x^4]=0$ if and only if the following conditions hold:
$$
k^2_1 k^4_4+k^2_5 k^4_1=0;\quad  k^2_2 k^4_4=0;\quad k^2_3
k^4_4=0;\quad k^2_7
k^4_1=0.
$$
So, the homologic elements with $x^0=0$ are of form
$$
x=y+z, \;\text{ where $y\in \Span(h^2_4, h^2_6, h^4_2, h^4_3, h^4_5,
h^4_6, h^4_7, h^6_1, h^6_2)$}
$$
and $z$ is one of the following:
$$
\renewcommand{\arraystretch}{1.4}
\begin{array}{l}
k^2_1h^2_1+ k^2_5h^2_5+ k^4_1h^4_1+ k^4_4h^4_4\;\text{ with
$k^2_1k^4_4+ k^2_5k^4_1=0$;}\\
k^2_2h^2_2+k^4_1h^4_1;\\
k^2_2h^2_2+k^2_7h^2_7;\\
k^2_3h^2_3+k^2_7h^2_7;\\
k^2_1h^2_1+ k^2_3h^2_3+k^4_1h^4_1;\\
k^2_5h^2_5+ k^2_7h^2_7+k^4_4h^4_4.
\end{array}
$$

b) $x^0\neq 0$, i.e., $k^0_1\neq 0$.  Then
since $[\mathbf{H}^0, \mathbf{H}^6]=0$, the condition $[x,x]=0$ holds
if and only if
$$
[x^0,x^2]=0;\quad [x^0,x^4]+[x^2,x^2]=0;\quad [x^2,x^4]=0.
$$
Using once more the multiplication table, we deduce that $x$ is homologic
if and only if
$$
x=y+z, \text{ where $y\in \Span(h^2_4,
h^2_6, h^4_2, h^4_6, h^4_7, h^6_1, h^6_2)$}
$$
and $z$ is only one of the following
(for $\alpha\neq 0$ and any $\beta,\ \gamma \in \Cee$):
$$
\renewcommand{\arraystretch}{1.4}
\begin{array}{l}
\alpha h^0_1+\beta
h^2_2+\gamma(h^2_3+540\frac{\beta}{\alpha}h^4_3);\\
\alpha
h^0_1+\beta h^2_2+\gamma(h^2_5-\frac{108\beta}{5\alpha}h^4_5).
\end{array}
$$

\ssec{2.5.2.  The homologic elements in $\fl_{\bcdot}$ for
$\fg=\fsl(3)$} Let $x\in\fl_{\bcdot}$ be a homologic element.  Then $x$ is
of the form $x=x^0+x^2$, cf. (7), and hence $x$ is homologic if and
only if
$$
[x^0,x^2]=0\quad \text{ and }\quad [x^2,x^2]=0.
$$

a) $x^0\neq 0$. One can see
from the multiplication table, that $[x^0,x^2]=0$  if and only if
$$
k^0_1 k^2_1= k^0_1 k^2_2= k^0_1 k^2_4=0,
$$
e.g., in this case
$k^2_1= k^2_2= k^2_4=0$. So, we need $x\in \Span(h^0_1, h^2_3,
h^2_5)$. One can also see from the table, that
$$
[h^0_1, h^2_3]=
[h^0_1, h^2_5]= [h^2_3, h^2_5]=0,
$$
i.e., any $x\in \Span(h^0_1,
h^2_3, h^2_5)$ is homologic. So, the final answer is:
$$
x=ah^0_1+bh^2_3+ch^2_5 \quad \text{for $a\neq 0$ and any $b, c\in\Cee$}.
$$

b) $x^0=0$.  Clearly, the elements that belong to either $L^{2;-2}$ or
$L^{1;-1}$ are homologic.  But there are also ``mixed'' elements $x^2$
satisfying $[x^2, x^2]=0$.  From the multiplication table we deduce
that $[x^2, x^2]=0$ if and only if
$$
k^2_1 k^2_5- k^2_2 k^2_3= k^2_2 k^2_5+ k^2_3
k^2_4=0.
$$
So, the answer in this case is:
$$
x=ah^2_1+ bh^2_2+ ch^2_3+ dh^2_4+ eh^2_5,\quad \text{where
$ae-bc=be+cd=0$}.
$$

\ssec{3.  A discussion} After a nice example of $\fl_{\bcdot}$ for
$\fg=\fsl(3)$ it was rather discouraging to discover the lack of
simple components in $\fl_{0}$ for all other algebras $\fg$.  Could it
be that the ``big" DG Lie superalgebra
$\fL_{\bcdot}=\Pi(C^{\bcdot+1})$ is more interesting than
$\fl_{\bcdot}$?  Let us consider two examples. 

\ssec{3.1.  The Nijenhuis bracket} Recall that the Nijenhuis bracket,
with which all similar examples started, is defined on the space of
sections of $\Omega^{\bcdot}\otimes_{\cF}\Vect$, where $\cF$, $\Vect$
and $\Omega^{\bcdot}$ are the sheaves of functions, vector fields and
differential forms on a given manifold.  Hereafter tensoring is
performed over $\cF$ and sheaves are replaced by modules and rings
of sections.

The {\it Nijenhuis bracket}
is defined, for any $\omega^k\in \Omega^k$, $\omega^l\in \Omega^l$ and
$\xi, \eta\in\Vect$, to be
$$
\renewcommand{\arraystretch}{1.4}
\begin{array}{l}
    \omega^k\otimes \xi, \omega^l\otimes \eta\mapsto
    (\omega^k\wedge \omega^l)\otimes [\xi, \eta]+\\
\left(\omega^k\wedge L_{\xi}(\omega^l)+(-1)^k   d\omega^k\wedge
\iota(\xi)(\omega^l)\right )\otimes \eta+\\
\left(-L_{\eta}(\omega^k)\wedge\omega^l+(-1)^l
\iota(\eta)(\omega^k)\wedge d\omega^l)\right
)\otimes \xi,
\end{array}
$$
where $\iota$ is the inner product and $L_{X}$ is the Lie derivative
with respect to the field $X$. The Nijenhuis
bracket has the following interpretation which implies its invariance:
the invariant operator $D: (\Omega^{k} \otimes_{\cF} \Vect(M),
\Omega^{\bcdot}) \tto \Omega^{\bcdot}$ given by the formula
$$
\renewcommand{\arraystretch}{1.4}
\begin{array}{l}
D(\omega^k\otimes \xi, \omega)=\\
d\left(\omega^k\wedge
\iota(\xi)(\omega)\right ) +(-1)^k\omega^k\wedge
\iota(\xi)(d\omega)=
d\omega^k\wedge
\iota(\xi)(\omega) +(-1)^k\omega^k\wedge
L_{\xi}(\omega)
\end{array}
$$
is, for a fixed $\omega^k\otimes \xi$, a superderivation of the
supercommutative superalgebra $\Omega^{\bcdot}$ and the Nijenhuis
bracket is just the supercommutator of these superderivations, see
\cite{G1}.  So we can identify $\Omega^{\bcdot} \otimes_{\cF}
\Vect(M)$, with the centralizer $C(d)$ of the exterior differential
considered as a vector field on $\hat M$, where $\hat M$ is the
supermanifold $(M, \Omega^{\bcdot}(M))$:
$$
C(d) = \{D\in\fvect(\hat M)\mid [D, d] = 0\}.
$$
(Here by $\fvect$ we denote the Lie superalgebra on the space of
sections of the sheaf $\Vect$.)

In contradistinction with a rich algebra of chains $\fL_{\bcdot}=C(d)$
we have $\fl_{\bcdot}=0$ (at least, locally). The above applies to
supermanifolds $M$ as well.

\ssec{3.2.  On formula $(1)$ and the codifferential $d$} Formula (1)
shows that $\fL_{\bcdot}$ is isomorphic to the simple Lie superalgebra
$\fvect(\dim \Pi(\fn))$ of polynomial vector fields on the superspace
$\Pi(\fn))^{*}$ {\it regardless} of multiplication in $\fn$.

The exterior differential in the cochain complex $C^{\bcdot}$ is just
a homological vector field of degree 1 in $\fvect(\dim \Pi(\fn)))$. 
The definition of the codifferential $d$ shows that if the bracket in
$\fn$ vanishes identically, as is the case for Hermitean symmetric
spaces $X=G/P$, the differential $d$ is also zero and
$\fl_{\bcdot}=\fvect(\dim \Pi(\fn))$.
 
In Example 3.1, $d=\sum\xi_{i}\pder{x_{i}}$, where $\xi_{i}= dx_{i}$. 
In this realization, $d$ is a degree 0 ``maximally 
nondegenerate'' (see \cite{V2}) vector field.  

If $\fn$ is a simple Lie algebra, the degree 1 vector field in
$\fvect(\dim \Pi(\fn)))$ is also ``maximally nondegenerate'' but,
unlike degree 0 fields, its explicit form is rather complicated,
except for $\fsl(3)$ and $\fgl(n)$, cf.  \cite{L2}.

Whatever the form of $d$, here is an implicit description of various
DG Lie superalgebras: in $\fvect(m|n)=\fder\Cee[X]$, where $X=(x_{1},
\dots, x_{m}, \xi_{1}, \dots, \xi_{n})$ and where $\deg~X_{i}=1$ for
all $i$, fix an element $d$ of degree 1.  Then the centralizer
$C(d)$ is a subalgebra in $\fvect(m|n)$ and $C(d)/\IM~d$ are DG Lie algebras. In 
sec. 2 we explicitely described $\fl_{\bcdot}=C(d)/\IM~d$ related with the 
nilpotent Lie algebra structure on $\Pi(\fvect(0|n)_{-1})$.

\end{document}